\documentclass[11pt,a4paper]{article}
\usepackage[T1]{fontenc}
\usepackage{lmodern}
\usepackage[margin=1in]{geometry}
\usepackage{amsmath,amssymb,amsthm}
\usepackage{microtype}
\usepackage{enumitem}
\usepackage{cite}
\usepackage[hidelinks]{hyperref}
\hypersetup{pdftitle={Understanding SAGA: Note VI on Variance Reduction},pdfauthor={Kamran Fareed and Alina Lodhi}}
\newtheorem{assumption}{Assumption}
\newtheorem{lemma}{Lemma}
\newtheorem{theorem}{Theorem}
\newcommand{\E}{\mathbb{E}}
\newcommand{\R}{\mathbb{R}}
\newcommand{\norm}[1]{\left\|#1\right\|}
\newcommand{\ip}[2]{\left\langle #1,#2\right\rangle}
\newcommand{\eps}{\varepsilon}
\setlist[itemize]{leftmargin=1.4em,itemsep=2pt,topsep=4pt,parsep=0pt}
\title{Variance Reduction Note 6: The Convergence of SAGA}
\author{Kamran Fareed \qquad Alina Lodhi\\
University of Agriculture Faisalabad\\
{\small kamrafareed96@gmail.com \qquad xieyunqi27@gmail.com}}
\date{}

\begin{document}
\maketitle
\vspace{-2em}
\begin{abstract}
SAGA corrects a sampled component gradient with a table of previously
computed gradients. Updating one entry at a time gives an unbiased
gradient estimate at the cost of one new component gradient per
iteration. This sixth note on variance reduction proves a linear
convergence guarantee for a smooth finite sum with convex components
and a strongly convex average. A two-component quadratic shows how the spread of the stored evaluation
points determines the variance, even at the minimizer. We track the
iterate and the table jointly to obtain a last-iterate bound, then
relate its dependence on the sample size to the rate at which individual
entries are replaced. For generalized linear losses the same gradient
table can be represented by one scalar per observation and a running mean.
\end{abstract}

\section{A control variate for each component}
Consider the finite-sum problem
\begin{equation}\label{eq:problem}
 \min_{x\in\R^d}F(x):=\frac1n\sum_{i=1}^n f_i(x).
\end{equation}
At a minimizer $x_\star$, the component gradients may point in
different directions even though their average vanishes. A control
variate should cancel these individual gradients as the iterates
approach $x_\star$. SVRG uses an exact gradient at a common snapshot
to do so \cite{svrg2013}. In Note~5, loopless SVRG refreshed that
snapshot with probability $1/n$, updating all the reference gradients
at once. SAGA, introduced by Defazio, Bach, and
Lacoste-Julien \cite{saga2014}, keeps a separate gradient $\alpha_i$
for each component. These gradients may have been computed at different
iterates. With $\bar\alpha=n^{-1}\sum_i\alpha_i$ and a uniform
random index $i$, the estimate at $x$ is
\begin{equation}\label{eq:estimator}
 v=\nabla f_i(x)-\alpha_i+\bar\alpha .
\end{equation}
For any fixed table, $\E_i[\alpha_i]=\bar\alpha$. Thus the correction
preserves the mean, however old the stored entries are. If the table
contains $\alpha_i=\nabla f_i(x_\star)$ for all $i$, then $v=0$ at
$x_\star$ for every choice of $i$.

The distinction from SAG \cite{sag2012} can be seen in one update.
Write $r_i=\nabla f_i(x)-\alpha_i$ for the correction to entry $i$.
After replacing the entry, the new table mean is
$\bar\alpha+r_i/n$. SAG uses this updated mean as its direction.
SAGA uses $\bar\alpha+r_i$ as its direction and reserves the factor
$1/n$ for maintaining the table mean. This distinction explains why
conditional unbiasedness is available in the proof below.

\section{Assumptions and the table update}
\begin{assumption}\label{ass:convex}
Each $f_i:\R^d\to\R$ is convex and differentiable, and its gradient
is $L$-Lipschitz for a common $L>0$. The average $F$ is
$\mu$-strongly convex for some $0<\mu\leq L$. In particular, for all
$x,y\in\R^d$,
\begin{equation}\label{eq:strong}
 F(y)\geq F(x)+\ip{\nabla F(x)}{y-x}+\frac{\mu}{2}\norm{y-x}^2.
\end{equation}
Write $x_\star$ for the unique minimizer of $F$,
$h(x)=F(x)-F(x_\star)$, and $\kappa=L/\mu$. Norms are Euclidean.
\end{assumption}
The average is also $L$-smooth. The components may have different
minimizers. Their convexity will be used to bound a component-gradient
difference by the objective gap after averaging over $i$.

Start at a deterministic $x_0$, compute
$\alpha_i^0=\nabla f_i(x_0)$ for each $i$, and form their mean
$\bar\alpha^0$. Fix a step size $\eta>0$. For $t=0,1,\ldots,T-1$,
perform the following steps:
\begin{enumerate}[leftmargin=1.6em,itemsep=2pt,topsep=4pt,parsep=0pt]
 \item Draw $i_t$ uniformly from $\{1,\ldots,n\}$, independently of
 the previous draws. Compute $g_t=\nabla f_{i_t}(x_t)$ and save
 $r_t=g_t-\alpha_{i_t}^t$.
 \item Using the mean before replacement, set
 \begin{equation}\label{eq:update}
 v_t=\bar\alpha^t+r_t,\qquad x_{t+1}=x_t-\eta v_t .
 \end{equation}
 \item Set $\alpha_{i_t}^{t+1}=g_t$ and $\alpha_j^{t+1}=\alpha_j^t$
 for $j\ne i_t$. Update $\bar\alpha^{t+1}=\bar\alpha^t+r_t/n$.
\end{enumerate}
Return $x_T$. The last update preserves the identity
$\bar\alpha^t=n^{-1}\sum_i\alpha_i^t$. Each iteration reads and
replaces one table entry and updates the mean in $O(d)$ arithmetic
operations in addition to the component-gradient computation. A dense
table uses $O(nd)$ storage. Initialization costs $n$ component-gradient
evaluations and the displayed implementation uses one further evaluation
per iteration. The table assumes that an individual component can be
revisited by its index.

\subsection*{A quadratic with equal component curvatures}
Consider two equally likely components on the real line,
\[
 f_+(x)=\tfrac12 ax^2+bx,\qquad
 f_-(x)=\tfrac12 ax^2-bx,\qquad a>0,\quad b\ne0 .
\]
Their average is $F(x)=ax^2/2$, with minimizer $x_\star=0$.
Both gradients change by $a(x-w)$ between any two points $w$ and $x$.
Thus the common-snapshot SVRG direction in Note~5 equals $ax$ for
either index, even when $w$ is far from $x$. For SAGA the two entries
can come from different points. Denote these points by $z_+,z_-$, so
that $\alpha_+=az_++b$ and $\alpha_-=az_--b$. Their mean is
$\bar\alpha=a\bar z$, where $\bar z=(z_++z_-)/2$. At the current point
$x$, substituting these entries into the SAGA update gives
\[
 v_i=ax+a(\bar z-z_i),\qquad
 \E_i[v_i]=ax,\qquad
 \E_i[(v_i-ax)^2]=\frac{a^2}{4}(z_+-z_-)^2 .
\]
The expectation here fixes $x$ and the table. The linear terms cancel,
although they give nonzero component gradients at the solution. The
conditional variance is determined by the separation of the two stored
points. Moving both stored points by the same amount leaves it unchanged.

For a short trajectory, take $a=b=1$, start at $x_0=4$, and use
$\eta=1/2$ for this example. Initially both entries are evaluated
at $4$. The first direction is $v_0=4$, so $x_1=2$. Updating either
entry at $x_0$ leaves the table as it was. At $t=1$ the direction is
again exact, giving $v_1=2$ and $x_2=1$. Suppose the index drawn at
$t=1$ was $+$. The stored points at the start of step $2$ are then
$(z_+,z_-)=(2,4)$. The next draw has two equally likely outcomes:
\[
 \begin{array}{c|c|c|c}
 i_2 & v_2 & x_3 & (z_+,z_-)\ \text{after the update}\\ \hline
 + & 2 & 0 & (1,4)\\
 - & 0 & 1 & (2,1)
 \end{array}
\]
At $x_2=1$ the two directions have mean $1$ and variance $1$.
On the first branch the iterate reaches the minimizer. The next
direction nevertheless takes the values $3/2$ and $-3/2$ with equal
probability, giving conditional variance $9/4$. On the second branch
the iterate stays at $1$, but the table changes. Its new separation is
one, so the next conditional variance is $1/4$. The first branch
shows why the distance to $x_\star$ alone cannot describe a SAGA step:
even when that distance is zero, the table can move the next iterate
away from the solution.

\subsection*{A scalar table for generalized linear models}
Following the scalar-storage observation in \cite{saga2014}, for a generalized linear model write $f_i(x)=\ell_i(a_i^\top x)$,
where $a_i\in\R^d$ is the feature vector of observation $i$ and the
label is included in the scalar loss $\ell_i$. The chain rule gives
$\nabla f_i(x)=\ell_i'(a_i^\top x)a_i$. Every gradient of component
$i$ is therefore a multiple of the same feature vector. To recover its
stored gradient $\alpha_i^t$, it suffices to keep the multiplier $s_i^t$
and read $a_i$ from the data. For example, with squared loss
$\ell_i(u)=(u-y_i)^2/2$, the stored number is the residual at the last
evaluation of observation $i$.

The initialization pass sets $s_i^0=\ell_i'(a_i^\top x_0)$ and accumulates
$\bar\alpha^0=n^{-1}\sum_i s_i^0a_i$. On drawing $i=i_t$, compute
$s=\ell_i'(a_i^\top x_t)$ and save the scalar difference
$\delta=s-s_i^t$. Using the old mean, form the direction and take the step
\[
 v_t=\bar\alpha^t+\delta a_i,\qquad x_{t+1}=x_t-\eta v_t.
\]
Now replace $s_i^t$ by $s$ and add $(\delta/n)a_i$ to the mean.
All other multipliers retain their values. Using the updated mean to
form $v_t$ would add an extra $(\delta/n)a_i$ to the direction.
The table and its mean require $n+d$ scalar entries in addition to the
data and iterate. If $a_i$ is sparse, the score and the change to the
mean use only its nonzero coordinates. The mean itself can be dense,
so an explicit update of $x_t$ can still touch all $d$ coordinates.

If the objective also contains a ridge penalty $\lambda\norm{x}^2/2$,
its gradient $\lambda x_t$ is known at every step. Keep the table for
the data losses and add $\lambda x_t$ to the displayed direction before
updating $x_t$. Its conditional expectation is then the gradient of the
regularized objective. The replacement of $s_i^t$ and the update to
$\bar\alpha^t$ still involve only the sampled loss.

\section{A short linear-convergence proof}
Let $\E_t$ denote conditional expectation given the history before
drawing $i_t$. Two quantities measure the accuracy of the fresh
gradient and the stored entries:
\begin{equation}\label{eq:errors}
 A(x)=\frac1n\sum_{i=1}^n\norm{\nabla f_i(x)-\nabla f_i(x_\star)}^2,
 \qquad
 H_t=\frac1n\sum_{i=1}^n\norm{\alpha_i^t-\nabla f_i(x_\star)}^2 .
\end{equation}
In particular $H_0=A(x_0)$. Both quantities refer to errors at the
solution and are used only in the analysis.

\begin{lemma}[A bound for fresh component gradients]\label{lem:growth}
Under Assumption~\ref{ass:convex}, $A(x)\leq2Lh(x)$ for every $x$.
\end{lemma}
\begin{proof}
Let $\psi_i(x)=f_i(x)-f_i(x_\star)
-\ip{\nabla f_i(x_\star)}{x-x_\star}$. Convexity makes $\psi_i$
nonnegative and its gradient is $L$-Lipschitz. Smoothness at
$y=x-\nabla\psi_i(x)/L$ gives
\[
 0\leq\psi_i(y)\leq\psi_i(x)-\frac{1}{2L}\norm{\nabla\psi_i(x)}^2 .
\]
Since $\nabla\psi_i(x)=\nabla f_i(x)-\nabla f_i(x_\star)$, averaging
this bound gives $A(x)\leq (2L/n)\sum_i\psi_i(x)=2Lh(x)$. The last
equality uses $\nabla F(x_\star)=0$.
\end{proof}

\begin{lemma}[Direction and table error]\label{lem:table}
The update \eqref{eq:update} satisfies
\begin{align}
 \E_t[v_t]&=\nabla F(x_t),&
 \E_t\norm{v_t}^2&\leq 4Lh(x_t)+2H_t,\label{eq:moment}\\
 \E_t H_{t+1}&=\left(1-\frac1n\right)H_t+\frac1n A(x_t).
 \label{eq:table}
\end{align}
\end{lemma}
\begin{proof}
The current table and $x_t$ are fixed under the conditioning.
Averaging \eqref{eq:estimator} proves the mean identity. Write
$a_i=\nabla f_i(x_t)-\nabla f_i(x_\star)$ and
$b_i=\alpha_i^t-\nabla f_i(x_\star)$. Since
$\bar b=n^{-1}\sum_i b_i=\bar\alpha^t$, the estimate is
$a_{i_t}-(b_{i_t}-\bar b)$. Therefore
\[
 \E_t\norm{v_t}^2\leq 2A(x_t)+2\E_t\norm{b_{i_t}-\bar b}^2
 \leq 2A(x_t)+2H_t .
\]
Lemma~\ref{lem:growth} gives \eqref{eq:moment}. Replacing entry $i_t$
changes $H_t$ by $(\norm{a_{i_t}}^2-\norm{b_{i_t}}^2)/n$.
Taking its conditional expectation proves \eqref{eq:table}.
\end{proof}

\begin{theorem}[Contraction at each iteration]\label{thm:linear}
Suppose Assumption~\ref{ass:convex} holds. Take $\eta=1/(6L)$ and define
\begin{equation}\label{eq:potential}
 \Phi_t=\norm{x_t-x_\star}^2+4n\eta^2 H_t,\qquad
 \gamma=\min\left\{\frac{\mu}{6L},\frac{1}{2n}\right\}.
\end{equation}
Then $\E_t\Phi_{t+1}\leq(1-\gamma)\Phi_t$. In particular, with
$C_0=L\Phi_0/2$, the last iterate after $T$ updates satisfies
\begin{equation}\label{eq:linear}
 \E h(x_T)\leq C_0(1-\gamma)^T .
\end{equation}
\end{theorem}
\begin{proof}
Put $D_t=\norm{x_t-x_\star}^2$. Strong convexity at $x_t$ with
$y=x_\star$ gives
$\ip{\nabla F(x_t)}{x_t-x_\star}\geq h(x_t)+\mu D_t/2$.
Expanding the squared distance in \eqref{eq:update} and applying
Lemma~\ref{lem:table} yields
\begin{align*}
 \E_t D_{t+1}
 &=D_t-2\eta\ip{\nabla F(x_t)}{x_t-x_\star}+\eta^2\E_t\norm{v_t}^2\\
 &\leq (1-\mu\eta)D_t+(-2\eta+4L\eta^2)h(x_t)+2\eta^2H_t .
\end{align*}
The term $2\eta^2H_t$ is paid for by the decrease of the old table
error. Multiply \eqref{eq:table} by $4n\eta^2$, bound
$A(x_t)\leq2Lh(x_t)$, and add the result to the distance inequality:
\[
 \E_t\Phi_{t+1}\leq(1-\mu\eta)D_t
 +4n\eta^2\left(1-\frac{1}{2n}\right)H_t
 +(-2\eta+12L\eta^2)h(x_t).
\]
At $\eta=1/(6L)$ the last coefficient vanishes. Both remaining
coefficients contract their respective terms by at least $\gamma$.
Induction using total expectation gives
$\E\Phi_T\leq (1-\gamma)^T\Phi_0$.
Finally $L$-smoothness at $x_\star$ and $\nabla F(x_\star)=0$
give $h(x_T)\leq LD_T/2\leq L\Phi_T/2$.
\end{proof}

\paragraph{Oracle complexity and accuracy.}
For an expected objective-gap target $0<\tau<C_0$, it suffices to take
\begin{equation}\label{eq:complexity}
 T=\left\lceil\max\{6\kappa,2n\}\log\frac{C_0}{\tau}\right\rceil,
 \qquad N=n+T .
\end{equation}
Indeed $(1-\gamma)^T\leq\exp(-\gamma T)$. If $\tau\geq C_0$,
$x_0$ already meets the target. The constant in \eqref{eq:linear}
includes the initial table. For the specified initialization,
\[
 C_0=\frac{L}{2}\norm{x_0-x_\star}^2+\frac{n}{18L}A(x_0)
 \leq\left(\kappa+\frac{n}{9}\right) h(x_0),
\]
using strong convexity and Lemma~\ref{lem:growth}. Thus a prescribed
iteration count may use an upper bound on $C_0$. For a relative target
$\E\Phi_T\leq\delta\Phi_0$, $0<\delta<1$, the count is
$n+\lceil\max\{6\kappa,2n\}\log(1/\delta)\rceil$. For an absolute
gradient target, the smoothness inequality at
$x-\nabla F(x)/L$ gives $\norm{\nabla F(x)}^2\leq2Lh(x)$.
Taking $\tau=\eps^2/(2L)$ ensures
$\E\norm{\nabla F(x_T)}^2\leq\eps^2$.

\subsection*{How long does an entry survive?}
An entry of the SAGA table has two histories worth distinguishing.
Its age records when the component was last sampled. Its contribution
to $H_t$ also depends on where the iterate was at that time. The
sampling rule determines the distribution of ages even before the
objective functions are specified. To see how the two histories enter
the proof, suppose $n\geq2$ and write $q=1-1/n$. For $n=1$, every draw
updates the sole entry.

Fix an index $j$. Conditional on the history at the start of any
step, the probability that the next $k$ draws all miss $j$ is $q^k$.
In particular, if $U_t$ counts the entries that have never been
selected in the first $t$ updates after initialization, then
\begin{equation}\label{eq:table-unvisited}
 \E U_t=nq^t,\qquad \Pr(U_t>0)\leq nq^t .
\end{equation}
The equality follows by adding the indicators of the $n$ individual
events. The inequality follows because $U_t\geq1$ whenever an entry
is still unvisited. Neither calculation requires independence between
the ages of different entries. After $n$ updates the expected
unvisited fraction is $(1-1/n)^n$, which tends to $e^{-1}$ as the
table grows. This count concerns selection: an entry selected at the
first step is visited even if its numerical value stays the same.

The same survival probability appears in the accumulated table error.
Taking expectations in \eqref{eq:table} and iterating the equality gives
\begin{equation}\label{eq:table-history}
 \E H_t=q^tH_0+\frac1n\sum_{s=0}^{t-1}q^{t-1-s}\E A(x_s).
\end{equation}
The first term accounts for entries whose values still come from
initialization. A gradient written at step $s$ survives to time $t$
if its index is missed in the following $t-1-s$ draws. This explains
the weight of $\E A(x_s)$ and why an update at an inaccurate iterate
can introduce fresh table error. The weights, including the initial
one, sum to one.

For comparison, let $\tau_{\rm all}$ be the first time all entries
have been selected. Once $k$ distinct indices have appeared, each
new draw finds an unseen one with probability $(n-k)/n$. Summing the
mean waiting times gives
$\E\tau_{\rm all}=\sum_{k=0}^{n-1}n/(n-k)=n\sum_{j=1}^n1/j$.
Full coverage therefore takes order $n\log n$ updates on average.
The contraction in Theorem~\ref{thm:linear} controls the average
squared error $H_t$ jointly with iterate distance. It applies from the
first update, including steps at which many entries still retain their
initial values.

\paragraph{Interpretation and scope.}
The two terms in $\gamma$ have different origins. Strong convexity
contracts the iterate distance on the scale of $\kappa$ updates. Each
table entry is refreshed with probability $1/n$ per update, so the
bound also depends on the time needed to replace stale gradients.
This proof uses uniform independent sampling, smooth convex components,
and a strongly convex average. The proximal and merely convex settings
in the original paper \cite{saga2014}, and extensions to arbitrary
sampling \cite{sagasampling2019}, require corresponding changes to the
update or its analysis.

\section{Recommended reading}

\paragraph{Variance-reduction methods.}
Basic guarantees \cite{ghadimi2013,gower2020,unifiedsgd2020}; gradient memories and control
variates \cite{sag2012,sdca2013,svrg2013,s2gd2013,saga2014,finito2014,proxsvrg2014}; sampling and acceleration \cite{catalyst2015,katyusha2016,svrgpp2016,katyushax2018,sagasampling2019,jacsketch2018,svrgsd2018,loopless2019,looplesssampling2021}; nonconvex and recursive
estimators \cite{reddi2016,allenzhu2016,natasha2017,sarah2017,isarah2018,scsg2017,spider2018,snvrg2018,spiderboost2019,proxsarah2020,zerosarah2021,page2021,nonconvexsampling2019}; momentum and adaptive schedules \cite{storm2019,hybrid2019,stormplus2021,superadam2021,metastorm2022,adaspider2022,jiang_adaptive2024,adasvrg2021,almosttune2020,adaptiveextra2022,alfpg2026};
second-order extensions \cite{ssrgd2019,natasha22018,svrc2019,svrcnewton2018}; lower bounds \cite{lower2023}.

\paragraph{Compositional optimization.}
Two-level  approximation \cite{scgd2014,asc2017,nasa2020,nearlyeasy2020}; finite-sum and recursive
estimators \cite{lian2017,compavg2018,civr2019,sarahcomp2019}; multilevel tracking \cite{nested2019,jiang_multilevel2022,levelindependent2020,jiang_revisiting2025};
coupled and partially observed blocks \cite{fccotheory2022,jiang_msvr2022,weakcomp2023,jiang_fcco2026,momentumcomp2025}; projection-free
methods \cite{vrfw2016,projectioncomp2022,jiang_projection2024,scgpp2019,jiang_projection2026}; statistical-model applications and federated
extensions \cite{jiang_models2023,fednest2022,jiang_models2026,localcomp2022}.

\paragraph{Sign-based optimization.}
Quantization, signSGD, majority voting, and error feedback \cite{qsgd2017,terngrad2017,signsgd2018,majority2019,errorfeedback2019};
variance-reduced sign updates and Lion \cite{jiang_sign2024,lion2023,jiang_lion2026}; compressed gradient
tracking \cite{diana2019,marina2021,jiang_distributed2026,ef212021,dasha2022}. Related distributed control-variate methods include
SCAFFOLD, FedPAGE, decentralized hybrid estimators, and ProxSkip \cite{scaffold2020,fedpage2021,hybriddecentral2021,proxskip2022}.

\paragraph{Minimax optimization.}
Finite-sum saddle-point methods, matrix games, and variance-reduced
extragradient \cite{saddle2016,matrixgames2019,svre2019,catalystminimax2020}; recursive, hybrid, adaptive, and decentralized
nonconvex minimax methods \cite{sreda2020,hybridminimax2020,adgda2023,dream2024}; a variational-inequality treatment \cite{vi2022}.



\begin{thebibliography}{100}
\small

\bibitem{svrg2013}
Rie Johnson and Tong Zhang.
\href{https://papers.nips.cc/paper_files/paper/2013/hash/ac1dd209cbcc5e5d1c6e28598e8cbbe8-Abstract.html}{Accelerating
stochastic gradient descent using predictive variance reduction}.
In \emph{NeurIPS}, volume 26, 2013.

\bibitem{saga2014}
Aaron Defazio, Francis Bach, and Simon Lacoste-Julien.
\href{https://arxiv.org/abs/1407.0202}{SAGA: A fast incremental gradient
method with support for non-strongly convex composite objectives}.
In \emph{NeurIPS}, volume 27, 2014.

\bibitem{sag2012}
Nicolas Le Roux, Mark Schmidt, and Francis Bach.
\href{https://arxiv.org/abs/1202.6258}{A stochastic gradient method with an
exponential convergence rate for finite training sets}.
In \emph{Advances in Neural Information Processing Systems (NeurIPS)},
volume 25, 2012.

\bibitem{sagasampling2019}
Xun Qian, Zheng Qu, and Peter Richtarik.
\href{https://proceedings.mlr.press/v97/qian19a.html}{SAGA with Arbitrary Sampling}.
In \emph{ICML}, PMLR 97:5190--5199, 2019.

\bibitem{ghadimi2013}
Saeed Ghadimi and Guanghui Lan.
\href{https://arxiv.org/abs/1309.5549}{Stochastic first- and zeroth-order
methods for nonconvex stochastic programming}.
\emph{SIAM Journal on Optimization}, 23(4):2341--2368, 2013.

\bibitem{gower2020}
Robert M. Gower, Mark Schmidt, Francis Bach, and Peter Richt\'{a}rik.
\href{https://arxiv.org/abs/2010.00892}{Variance-reduced methods for machine
learning}.
\emph{Proceedings of the IEEE}, 108(11):1968--1983, 2020.

\bibitem{unifiedsgd2020}
Eduard Gorbunov, Filip Hanzely, and Peter Richtarik.
\href{https://proceedings.mlr.press/v108/gorbunov20a.html}{A Unified Theory of SGD: Variance Reduction, Sampling, Quantization and Coordinate Descent}.
In \emph{AISTATS}, PMLR 108:680--690, 2020.

\bibitem{sdca2013}
Shai Shalev-Shwartz and Tong Zhang.
\href{https://www.jmlr.org/papers/v14/shalev-shwartz13a.html}{Stochastic Dual Coordinate Ascent Methods for Regularized Loss Minimization}.
\emph{Journal of Machine Learning Research}, 14(16):567--599, 2013.

\bibitem{s2gd2013}
Jakub Konecny and Peter Richtarik.
\href{https://arxiv.org/abs/1312.1666}{Semi-Stochastic Gradient Descent Methods}.
\emph{arXiv preprint arXiv:1312.1666}, 2013.

\bibitem{finito2014}
Aaron J. Defazio, Tiberio S. Caetano, and Justin Domke.
\href{https://arxiv.org/abs/1407.2710}{Finito: A Faster, Permutable Incremental Gradient Method for Big Data Problems}.
\emph{arXiv preprint arXiv:1407.2710}, 2014.

\bibitem{proxsvrg2014}
Lin Xiao and Tong Zhang.
\href{https://arxiv.org/abs/1403.4699}{A proximal stochastic gradient
method with progressive variance reduction}.
\emph{SIAM Journal on Optimization}, 24(4):2057--2075, 2014.

\bibitem{catalyst2015}
Hongzhou Lin, Julien Mairal, and Zaid Harchaoui.
\href{https://papers.nips.cc/paper/2015/hash/c164bbc9d6c72a52c599bbb43d8db8e1-Abstract.html}{A Universal Catalyst for First-Order Optimization}.
In \emph{NeurIPS}, volume 28, 2015.

\bibitem{katyusha2016}
Zeyuan Allen-Zhu.
\href{https://arxiv.org/abs/1603.05953}{Katyusha: The first direct
acceleration of stochastic gradient methods}.
\emph{arXiv preprint arXiv:1603.05953}, 2016.

\bibitem{svrgpp2016}
Zeyuan Allen-Zhu and Yang Yuan.
\href{https://proceedings.mlr.press/v48/allen-zhub16.html}{Improved SVRG for Non-Strongly-Convex or Sum-of-Non-Convex Objectives}.
In \emph{ICML}, PMLR 48:1080--1089, 2016.

\bibitem{katyushax2018}
Zeyuan Allen-Zhu.
\href{https://proceedings.mlr.press/v80/allen-zhu18a.html}{Katyusha X: Simple Momentum Method for Stochastic Sum-of-Nonconvex Optimization}.
In \emph{ICML}, PMLR 80:179--185, 2018.

\bibitem{jacsketch2018}
Robert M. Gower, Peter Richtarik, and Francis Bach.
\href{https://arxiv.org/abs/1805.02632}{Stochastic Quasi-Gradient Methods: Variance Reduction via Jacobian Sketching}.
\emph{arXiv preprint arXiv:1805.02632}, 2018.

\bibitem{svrgsd2018}
Fanhua Shang, Yuanyuan Liu, Kaiwen Zhou, James Cheng, Kelvin Kai Wing Ng, and Yuichi Yoshida.
\href{https://proceedings.mlr.press/v84/shang18a.html}{Guaranteed Sufficient Decrease for Stochastic Variance Reduced Gradient Optimization}.
In \emph{AISTATS}, PMLR 84:1027--1036, 2018.

\bibitem{loopless2019}
Dmitry Kovalev, Samuel Horvath, and Peter Richtarik.
\href{https://arxiv.org/abs/1901.08689}{Don't Jump Through Hoops and Remove Those Loops: SVRG and Katyusha are Better Without the Outer Loop}.
\emph{arXiv preprint arXiv:1901.08689}, 2019.

\bibitem{looplesssampling2021}
Xun Qian, Zheng Qu, and Peter Richtarik.
\href{https://jmlr.org/papers/v22/20-156.html}{L-SVRG and L-Katyusha with Arbitrary Sampling}.
\emph{Journal of Machine Learning Research}, 22(112), 2021.

\bibitem{reddi2016}
Sashank J. Reddi, Ahmed Hefny, Suvrit Sra, Barnabas Poczos, and Alex Smola.
\href{https://proceedings.mlr.press/v48/reddi16.html}{Stochastic variance
reduction for nonconvex optimization}.
In \emph{International Conference on Machine Learning (ICML)},
PMLR 48:314--323, 2016.

\bibitem{allenzhu2016}
Zeyuan Allen-Zhu and Elad Hazan.
\href{https://proceedings.mlr.press/v48/allen-zhua16.html}{Variance reduction
for faster non-convex optimization}.
In \emph{ICML}, PMLR 48:699--707, 2016.

\bibitem{natasha2017}
Zeyuan Allen-Zhu.
\href{https://proceedings.mlr.press/v70/allen-zhu17a.html}{Natasha: Faster Non-Convex Stochastic Optimization via Strongly Non-Convex Parameter}.
In \emph{ICML}, PMLR 70:89--97, 2017.

\bibitem{sarah2017}
Lam M. Nguyen, Jie Liu, Katya Scheinberg, and Martin Tak\'{a}\v{c}.
\href{https://arxiv.org/abs/1703.00102}{SARAH: A novel method for machine
learning problems using stochastic recursive gradient}.
In \emph{ICML}, 2017.

\bibitem{isarah2018}
Lam M. Nguyen, Katya Scheinberg, and Martin Takac.
\href{https://arxiv.org/abs/1811.10105}{Inexact SARAH Algorithm for Stochastic Optimization}.
\emph{arXiv preprint arXiv:1811.10105}, 2018.

\bibitem{scsg2017}
Lihua Lei, Cheng Ju, Jianbo Chen, and Michael I. Jordan.
\href{https://arxiv.org/abs/1706.09156}{Non-convex finite-sum optimization
via SCSG methods}.
In \emph{NeurIPS}, volume 30, 2017.

\bibitem{spider2018}
Cong Fang, Chris Junchi Li, Zhouchen Lin, and Tong Zhang.
\href{https://arxiv.org/abs/1807.01695}{SPIDER: Near-optimal non-convex
optimization via stochastic path-integrated differential estimator}.
In \emph{NeurIPS}, volume 31, 2018.

\bibitem{snvrg2018}
Dongruo Zhou, Pan Xu, and Quanquan Gu.
\href{https://arxiv.org/abs/1806.07811}{Stochastic nested variance reduction
for nonconvex optimization}.
In \emph{NeurIPS}, volume 31, 2018.

\bibitem{spiderboost2019}
Zhe Wang, Kaiyi Ji, Yi Zhou, Yingbin Liang, and Vahid Tarokh.
\href{https://arxiv.org/abs/1810.10690}{SpiderBoost and momentum: Faster
stochastic variance reduction algorithms}.
In \emph{NeurIPS}, volume 32, 2019.

\bibitem{proxsarah2020}
Nhan H. Pham, Lam M. Nguyen, Dzung T. Phan, and Quoc Tran-Dinh.
\href{https://www.jmlr.org/papers/v21/19-248.html}{ProxSARAH: An Efficient Algorithmic Framework for Stochastic Composite Nonconvex Optimization}.
\emph{Journal of Machine Learning Research}, 21(110):1--48, 2020.

\bibitem{zerosarah2021}
Zhize Li, Slavomir Hanzely, and Peter Richtarik.
\href{https://arxiv.org/abs/2103.01447}{ZeroSARAH: Efficient Nonconvex Finite-Sum Optimization with Zero Full Gradient Computation}.
\emph{arXiv preprint arXiv:2103.01447}, 2021.

\bibitem{page2021}
Zhize Li, Hongyan Bao, Xiangliang Zhang, and Peter Richt\'{a}rik.
\href{https://proceedings.mlr.press/v139/li21a.html}{PAGE: A simple and
optimal probabilistic gradient estimator for nonconvex optimization}.
In \emph{ICML}, 2021.

\bibitem{nonconvexsampling2019}
Samuel Horvath and Peter Richtarik.
\href{https://proceedings.mlr.press/v97/horvath19a.html}{Nonconvex Variance Reduced Optimization with Arbitrary Sampling}.
In \emph{ICML}, PMLR 97:2781--2789, 2019.

\bibitem{storm2019}
Ashok Cutkosky and Francesco Orabona.
\href{https://arxiv.org/abs/1905.10018}{Momentum-based variance reduction
in non-convex SGD}.
In \emph{NeurIPS}, volume 32, 2019.

\bibitem{hybrid2019}
Quoc Tran-Dinh, Nhan H. Pham, Dzung T. Phan, and Lam M. Nguyen.
\href{https://arxiv.org/abs/1907.03793}{A hybrid stochastic optimization
framework for stochastic composite nonconvex optimization}.
\emph{arXiv preprint arXiv:1907.03793}, 2019.

\bibitem{stormplus2021}
Kfir Y. Levy, Ali Kavis, and Volkan Cevher.
\href{https://proceedings.neurips.cc/paper/2021/hash/ac10ff1941c540cd87c107330996f4f6-Abstract.html}{STORM+:
Fully adaptive SGD with recursive momentum for nonconvex optimization}.
In \emph{NeurIPS}, volume 34, 2021.

\bibitem{superadam2021}
Feihu Huang, Junyi Li, and Heng Huang.
\href{https://proceedings.neurips.cc/paper/2021/hash/4be5a36cbaca8ab9d2066debfe4e65c1-Abstract.html}{SUPER-ADAM: Faster and Universal Framework of Adaptive Gradients}.
In \emph{NeurIPS}, volume 34, 2021.

\bibitem{metastorm2022}
Zijian Liu, Ta Duy Nguyen, Thien Hang Nguyen, Alina Ene, and Huy L. Nguyen.
\href{https://arxiv.org/abs/2209.14853}{META-STORM: Generalized fully-adaptive
variance reduced SGD for unbounded functions}.
\emph{arXiv preprint arXiv:2209.14853}, 2022.

\bibitem{adaspider2022}
Ali Kavis, Stratis Skoulakis, Kimon Antonakopoulos, Leello Tadesse Dadi, and Volkan Cevher.
\href{https://proceedings.neurips.cc/paper_files/paper/2022/hash/94f625dcdec313cd432d65f96fcc51c8-Abstract-Conference.html}{Adaptive Stochastic Variance Reduction for Non-convex Finite-Sum Minimization}.
In \emph{NeurIPS}, volume 35, 2022.

\bibitem{jiang_adaptive2024}
Wei Jiang, Sifan Yang, Yibo Wang, and Lijun Zhang.
\href{https://proceedings.neurips.cc/paper_files/paper/2024/hash/272efd3a6091ceefcbc79f1f3a6fdba4-Abstract-Conference.html}{Adaptive
variance reduction for stochastic optimization under weaker assumptions}.
In \emph{NeurIPS}, volume 37, pages 22047--22080, 2024.

\bibitem{adasvrg2021}
Benjamin Dubois-Taine, Sharan Vaswani, Reza Babanezhad, Mark Schmidt, and Simon Lacoste-Julien.
\href{https://arxiv.org/abs/2102.09645}{SVRG Meets AdaGrad: Painless Variance Reduction}.
\emph{arXiv preprint arXiv:2102.09645}, 2021.

\bibitem{almosttune2020}
Bingcong Li, Lingda Wang, and Georgios B. Giannakis.
\href{https://proceedings.mlr.press/v119/li20n.html}{Almost Tune-Free Variance Reduction}.
In \emph{ICML}, PMLR 119:5969--5978, 2020.

\bibitem{adaptiveextra2022}
Zijian Liu, Ta Duy Nguyen, Alina Ene, and Huy L. Nguyen.
\href{https://proceedings.mlr.press/v162/liu22o.html}{Adaptive Accelerated (Extra-)Gradient Methods with Variance Reduction}.
In \emph{ICML}, PMLR 162:13947--13994, 2022.

\bibitem{alfpg2026}
Ganzhao Yuan.
\href{https://arxiv.org/abs/2502.21099v4}{Adaptive Lipschitz-Free Proximal Gradient Methods for Stochastic Composite Nonconvex Optimization}.
\emph{arXiv preprint arXiv:2502.21099v4}, 2026.

\bibitem{ssrgd2019}
Zhize Li.
\href{https://papers.nips.cc/paper_files/paper/2019/hash/addfa9b7e234254d26e9c7f2af1005cb-Abstract.html}{SSRGD: Simple Stochastic Recursive Gradient Descent for Escaping Saddle Points}.
In \emph{NeurIPS}, volume 32, 2019.

\bibitem{natasha22018}
Zeyuan Allen-Zhu.
\href{https://proceedings.neurips.cc/paper/2018/hash/79a49b3e3762632813f9e35f4ba53d6c-Abstract.html}{Natasha 2: Faster Non-Convex Optimization Than SGD}.
In \emph{NeurIPS}, volume 31, 2018.

\bibitem{svrc2019}
Zhe Wang, Yi Zhou, Yingbin Liang, and Guanghui Lan.
\href{https://proceedings.mlr.press/v89/wang19d.html}{Stochastic Variance-Reduced Cubic Regularization for Nonconvex Optimization}.
In \emph{AISTATS}, PMLR 89:2731--2740, 2019.

\bibitem{svrcnewton2018}
Dongruo Zhou, Pan Xu, and Quanquan Gu.
\href{https://proceedings.mlr.press/v80/zhou18d.html}{Stochastic Variance-Reduced Cubic Regularized Newton Methods}.
In \emph{ICML}, PMLR 80:5990--5999, 2018.

\bibitem{lower2023}
Yossi Arjevani, Yair Carmon, John C. Duchi, Dylan J. Foster,
Nathan Srebro, and Blake Woodworth.
\href{https://doi.org/10.1007/s10107-022-01822-7}{Lower bounds for non-convex
stochastic optimization}.
\emph{Mathematical Programming}, 199:165--214, 2023.

\bibitem{scgd2014}
Mengdi Wang, Ethan X. Fang, and Han Liu.
\href{https://arxiv.org/abs/1411.3803}{Stochastic compositional gradient
descent: Algorithms for minimizing compositions of expected-value functions}.
\emph{arXiv preprint arXiv:1411.3803}, 2014.

\bibitem{asc2017}
Mengdi Wang, Ji Liu, and Ethan X. Fang.
\href{https://jmlr.org/papers/v18/16-504.html}{Accelerating stochastic
composition optimization}.
\emph{Journal of Machine Learning Research}, 18(105):1--23, 2017.

\bibitem{nasa2020}
Saeed Ghadimi, Andrzej Ruszczynski, and Mengdi Wang.
\href{https://doi.org/10.1137/18M1230542}{A Single Timescale Stochastic Approximation Method for Nested Stochastic Optimization}.
\emph{SIAM Journal on Optimization}, 2020.

\bibitem{nearlyeasy2020}
Tianyi Chen, Yuejiao Sun, and Wotao Yin.
\href{https://arxiv.org/abs/2008.10847}{Solving Stochastic Compositional Optimization is Nearly as Easy as Solving Stochastic Optimization}.
\emph{arXiv preprint arXiv:2008.10847}, 2020.

\bibitem{lian2017}
Xiangru Lian, Mengdi Wang, and Ji Liu.
\href{https://proceedings.mlr.press/v54/lian17a.html}{Finite-sum composition
optimization via variance reduced gradient descent}.
In \emph{International Conference on Artificial Intelligence and Statistics
(AISTATS)}, PMLR 54:1159--1167, 2017.

\bibitem{compavg2018}
Tsung-Yu Hsieh, Yasser El-Manzalawy, Yiwei Sun, and Vasant Honavar.
\href{https://arxiv.org/abs/1809.01225}{Compositional Stochastic Average Gradient for Machine Learning and Related Applications}.
\emph{arXiv preprint arXiv:1809.01225}, 2018.

\bibitem{civr2019}
Junyu Zhang and Lin Xiao.
\href{https://proceedings.mlr.press/v97/zhang19n.html}{A composite randomized
incremental gradient method}.
In \emph{ICML}, PMLR 97:7454--7462, 2019.

\bibitem{sarahcomp2019}
Wenqing Hu, Chris Junchi Li, Xiangru Lian, Ji Liu, and Huizhuo Yuan.
\href{https://papers.neurips.cc/paper_files/paper/2019/hash/21ce689121e39821d07d04faab328370-Abstract.html}{Efficient Smooth Non-Convex Stochastic Compositional Optimization via Stochastic Recursive Gradient Descent}.
In \emph{NeurIPS}, volume 32, 2019.

\bibitem{nested2019}
Junyu Zhang and Lin Xiao.
\href{https://arxiv.org/abs/1908.11468}{Multi-level composite stochastic
optimization via nested variance reduction}.
\emph{arXiv preprint arXiv:1908.11468}, 2019.

\bibitem{jiang_multilevel2022}
Wei Jiang, Bokun Wang, Yibo Wang, Lijun Zhang, and Tianbao Yang.
\href{https://proceedings.mlr.press/v162/jiang22c.html}{Optimal algorithms
for stochastic multi-level compositional optimization}.
In \emph{ICML}, PMLR 162:10195--10216, 2022.

\bibitem{levelindependent2020}
Krishnakumar Balasubramanian, Saeed Ghadimi, and Anthony Nguyen.
\href{https://arxiv.org/abs/2008.10526}{Stochastic Multi-level Composition Optimization Algorithms with Level-Independent Convergence Rates}.
\emph{arXiv preprint arXiv:2008.10526}, 2020.

\bibitem{jiang_revisiting2025}
Wei Jiang, Sifan Yang, Yibo Wang, Tianbao Yang, and Lijun Zhang.
\href{https://doi.org/10.1109/TPAMI.2025.3552197}{Revisiting stochastic
multi-level compositional optimization}.
\emph{IEEE Transactions on Pattern Analysis and Machine Intelligence},
47(7):5613--5624, 2025.

\bibitem{fccotheory2022}
Bokun Wang and Tianbao Yang.
\href{https://proceedings.mlr.press/v162/wang22ak.html}{Finite-Sum Coupled Compositional Stochastic Optimization: Theory and Applications}.
In \emph{ICML}, PMLR 162:23292--23317, 2022.

\bibitem{jiang_msvr2022}
Wei Jiang, Gang Li, Yibo Wang, Lijun Zhang, and Tianbao Yang.
\href{https://papers.nips.cc/paper_files/paper/2022/hash/d13ee73683fd5567e5c07634a25cd7b8-Abstract-Conference.html}{Multi-block-single-probe
variance reduced estimator for coupled compositional optimization}.
In \emph{NeurIPS}, volume 35, pages 32499--32511, 2022.

\bibitem{weakcomp2023}
Quanqi Hu, Dixian Zhu, and Tianbao Yang.
\href{https://papers.nips.cc/paper_files/paper/2023/hash/1160792eab11de2bbaf9e71fce191e8c-Abstract-Conference.html}{Non-Smooth Weakly-Convex Finite-sum Coupled Compositional Optimization}.
In \emph{NeurIPS}, volume 36, 2023.

\bibitem{jiang_fcco2026}
Wei Jiang, Sifan Yang, Yibo Wang, Lijun Zhang, and Zechao Li.
\href{https://arxiv.org/abs/2609.15723}{Solving finite-sum coupled
compositional optimization via multi-block-single-probe estimator}.
\emph{arXiv preprint arXiv:2609.15723}, 2026.

\bibitem{momentumcomp2025}
Xingyu Chen, Bokun Wang, Ming Yang, Qihang Lin, and Tianbao Yang.
\href{https://papers.nips.cc/paper_files/paper/2025/hash/714467868756e0e0614e99de02ac0032-Abstract-Conference.html}{Stochastic Momentum Methods for Non-smooth Non-Convex Finite-Sum Coupled Compositional Optimization}.
In \emph{NeurIPS}, volume 38, 2025.

\bibitem{vrfw2016}
Elad Hazan and Haipeng Luo.
\href{https://proceedings.mlr.press/v48/hazana16.html}{Variance-Reduced and Projection-Free Stochastic Optimization}.
In \emph{ICML}, PMLR 48:1263--1271, 2016.

\bibitem{projectioncomp2022}
Tesi Xiao, Krishnakumar Balasubramanian, and Saeed Ghadimi.
\href{https://proceedings.neurips.cc/paper_files/paper/2022/hash/7e16384b94a1c7e4462a70bb8fb93ca9-Abstract-Conference.html}{A Projection-free Algorithm for Constrained Stochastic Multi-level Composition Optimization}.
In \emph{NeurIPS}, volume 35, 2022.

\bibitem{jiang_projection2024}
Wei Jiang, Sifan Yang, Wenhao Yang, Yibo Wang, Yuanyu Wan, and Lijun Zhang.
\href{https://proceedings.mlr.press/v235/jiang24g.html}{Projection-free
variance reduction methods for stochastic constrained multi-level
compositional optimization}.
In \emph{ICML}, PMLR 235:21962--21987, 2024.

\bibitem{scgpp2019}
Hamed Hassani, Amin Karbasi, Aryan Mokhtari, and Zebang Shen.
\href{https://arxiv.org/abs/1902.06992}{Stochastic Conditional Gradient++}.
\emph{arXiv preprint arXiv:1902.06992}, 2019.

\bibitem{jiang_projection2026}
Wei Jiang, Sifan Yang, Wenhao Yang, Yibo Wang, Yuanyu Wan,
Zechao Li, and Lijun Zhang.
\href{https://arxiv.org/abs/2609.15679}{Projection-free multi-level
algorithms for stochastic constrained compositional optimization}.
\emph{arXiv preprint arXiv:2609.15679}, 2026.

\bibitem{jiang_models2023}
Wei Jiang, Jiayu Qin, Lingyu Wu, Changyou Chen, Tianbao Yang, and Lijun Zhang.
\href{https://proceedings.mlr.press/v202/jiang23g.html}{Learning unnormalized
statistical models via compositional optimization}.
In \emph{ICML}, PMLR 202:15105--15124, 2023.

\bibitem{fednest2022}
Davoud Ataee Tarzanagh, Mingchen Li, Christos Thrampoulidis, and Samet Oymak.
\href{https://proceedings.mlr.press/v162/tarzanagh22a.html}{FedNest: Federated Bilevel, Minimax, and Compositional Optimization}.
In \emph{ICML}, PMLR 162:21146--21179, 2022.

\bibitem{jiang_models2026}
Wei Jiang, Jiayu Qin, Lingyu Wu, Changyou Chen, Tianbao Yang, and Lijun Zhang.
\href{https://jiangw1998.github.io/pdf/TPAMI-2026-Jiang.pdf}{Optimizing
unnormalized statistical models through compositional optimization}.
\emph{IEEE Transactions on Pattern Analysis and Machine Intelligence},
48(2):1949--1960, 2026.

\bibitem{localcomp2022}
Hongchang Gao, Junyi Li, and Heng Huang.
\href{https://proceedings.mlr.press/v162/gao22c.html}{On the Convergence of Local Stochastic Compositional Gradient Descent with Momentum}.
In \emph{ICML}, PMLR 162:7017--7035, 2022.

\bibitem{qsgd2017}
Dan Alistarh, Demjan Grubic, Jerry Li, Ryota Tomioka, and Milan Vojnovic.
\href{https://papers.nips.cc/paper/2017/hash/6c340f25839e6acdc73414517203f5f0-Abstract.html}{QSGD: Communication-Efficient SGD via Gradient Quantization and Encoding}.
In \emph{NeurIPS}, volume 30, 2017.

\bibitem{terngrad2017}
Wei Wen, Cong Xu, Feng Yan, Chunpeng Wu, Yandan Wang, Yiran Chen, and Hai Li.
\href{https://arxiv.org/abs/1705.07878}{TernGrad: Ternary Gradients to Reduce Communication in Distributed Deep Learning}.
In \emph{NeurIPS}, volume 30, 2017.

\bibitem{signsgd2018}
Jeremy Bernstein, Yu-Xiang Wang, Kamyar Azizzadenesheli, and Anima Anandkumar.
\href{https://proceedings.mlr.press/v80/bernstein18a.html}{signSGD:
Compressed optimisation for non-convex problems}.
In \emph{ICML}, PMLR 80:560--569, 2018.

\bibitem{majority2019}
Jeremy Bernstein, Jiawei Zhao, Kamyar Azizzadenesheli, and Anima Anandkumar.
\href{https://arxiv.org/abs/1810.05291}{signSGD with majority vote is
communication efficient and fault tolerant}.
In \emph{International Conference on Learning Representations (ICLR)}, 2019.

\bibitem{errorfeedback2019}
Sai Praneeth Karimireddy, Quentin Rebjock, Sebastian U. Stich, and Martin Jaggi.
\href{https://proceedings.mlr.press/v97/karimireddy19a.html}{Error feedback
fixes SignSGD and other gradient compression schemes}.
In \emph{ICML}, PMLR 97:3252--3261, 2019.

\bibitem{jiang_sign2024}
Wei Jiang, Sifan Yang, Wenhao Yang, and Lijun Zhang.
\href{https://papers.nips.cc/paper_files/paper/2024/hash/3b7db05c0e383518d789b6e93131f1f0-Abstract-Conference.html}{Efficient
sign-based optimization: Accelerating convergence via variance reduction}.
In \emph{NeurIPS}, volume 37, pages 33891--33932, 2024.

\bibitem{lion2023}
Xiangning Chen, Chen Liang, Da Huang, Esteban Real, Kaiyuan Wang,
Hieu Pham, Xuanyi Dong, Thang Luong, Cho-Jui Hsieh, Yifeng Lu, and Quoc V. Le.
\href{https://proceedings.neurips.cc/paper_files/paper/2023/hash/9a39b4925e35cf447ccba8757137d84f-Abstract.html}{Symbolic
discovery of optimization algorithms}.
In \emph{NeurIPS}, volume 36, 2023.

\bibitem{jiang_lion2026}
Wei Jiang, Mao Xu, Wenhao Yang, Yibo Wang, Zechao Li, and Lijun Zhang.
\href{https://proceedings.mlr.press/v306/jiang26ax.html}{Convergence analysis
of the Lion optimizer in centralized and distributed settings}.
In \emph{ICML}, PMLR 306:53308--53333, 2026.

\bibitem{diana2019}
Konstantin Mishchenko, Eduard Gorbunov, Martin Takac, and Peter Richtarik.
\href{https://arxiv.org/abs/1901.09269}{Distributed Learning with Compressed Gradient Differences}.
\emph{arXiv preprint arXiv:1901.09269}, 2019.

\bibitem{marina2021}
Eduard Gorbunov, Konstantin P. Burlachenko, Zhize Li, and Peter Richtarik.
\href{https://proceedings.mlr.press/v139/gorbunov21a.html}{MARINA: Faster Non-Convex Distributed Learning with Compression}.
In \emph{ICML}, PMLR 139:3788--3798, 2021.

\bibitem{jiang_distributed2026}
Wei Jiang, Zechao Li, and Lijun Zhang.
\href{https://arxiv.org/abs/2609.18656}{Revisiting distributed sign-based
variance reduction}.
\emph{arXiv preprint arXiv:2609.18656}, 2026.

\bibitem{ef212021}
Peter Richtarik, Igor Sokolov, and Ilyas Fatkhullin.
\href{https://papers.nips.cc/paper_files/paper/2021/hash/231141b34c82aa95e48810a9d1b33a79-Abstract.html}{EF21: A New, Simpler, Theoretically Better, and Practically Faster Error Feedback}.
In \emph{NeurIPS}, volume 34, 2021.

\bibitem{dasha2022}
Alexander Tyurin and Peter Richtarik.
\href{https://arxiv.org/abs/2202.01268}{DASHA: Distributed Nonconvex Optimization with Communication Compression, Optimal Oracle Complexity, and No Client Synchronization}.
\emph{arXiv preprint arXiv:2202.01268}, 2022.

\bibitem{scaffold2020}
Sai Praneeth Karimireddy, Satyen Kale, Mehryar Mohri, Sashank Reddi, Sebastian Stich, and Ananda Theertha Suresh.
\href{https://proceedings.mlr.press/v119/karimireddy20a.html}{SCAFFOLD: Stochastic Controlled Averaging for Federated Learning}.
In \emph{ICML}, PMLR 119:5132--5143, 2020.

\bibitem{fedpage2021}
Haoyu Zhao, Zhize Li, and Peter Richtarik.
\href{https://arxiv.org/abs/2108.04755}{FedPAGE: A Fast Local Stochastic Gradient Method for Communication-Efficient Federated Learning}.
\emph{arXiv preprint arXiv:2108.04755}, 2021.

\bibitem{hybriddecentral2021}
Ran Xin, Usman A. Khan, and Soummya Kar.
\href{https://proceedings.mlr.press/v139/xin21a.html}{A Hybrid Variance-Reduced Method for Decentralized Stochastic Non-Convex Optimization}.
In \emph{ICML}, PMLR 139:11459--11469, 2021.

\bibitem{proxskip2022}
Konstantin Mishchenko, Grigory Malinovsky, Sebastian Stich, and Peter Richtarik.
\href{https://proceedings.mlr.press/v162/mishchenko22b.html}{ProxSkip: Yes! Local Gradient Steps Provably Lead to Communication Acceleration! Finally!}.
In \emph{ICML}, PMLR 162:15750--15769, 2022.

\bibitem{saddle2016}
P. Balamurugan and Francis Bach.
\href{https://arxiv.org/abs/1605.06398}{Stochastic variance reduction
methods for saddle-point problems}.
In \emph{NeurIPS}, volume 29, 2016.

\bibitem{matrixgames2019}
Yair Carmon, Yujia Jin, Aaron Sidford, and Kevin Tian.
\href{https://papers.nips.cc/paper_files/paper/2019/hash/6c442e0e996fa84f344a14927703a8c1-Abstract.html}{Variance Reduction for Matrix Games}.
In \emph{NeurIPS}, volume 32, 2019.

\bibitem{svre2019}
Tatjana Chavdarova, Gauthier Gidel, Francois Fleuret, and Simon Lacoste-Julien.
\href{https://papers.nips.cc/paper_files/paper/2019/hash/58a2fc6ed39fd083f55d4182bf88826d-Abstract.html}{Reducing Noise in GAN Training with Variance Reduced Extragradient}.
In \emph{NeurIPS}, volume 32, 2019.

\bibitem{catalystminimax2020}
Junchi Yang, Siqi Zhang, Negar Kiyavash, and Niao He.
\href{https://papers.nips.cc/paper_files/paper/2020/hash/3db54f5573cd617a0112d35dd1e6b1ef-Abstract.html}{A Catalyst Framework for Minimax Optimization}.
In \emph{NeurIPS}, volume 33, 2020.

\bibitem{sreda2020}
Luo Luo, Haishan Ye, Zhichao Huang, and Tong Zhang.
\href{https://papers.nips.cc/paper/2020/hash/ecb47fbb07a752413640f82a945530f8-Abstract.html}{Stochastic
recursive gradient descent ascent for stochastic nonconvex-strongly-concave
minimax problems}.
In \emph{NeurIPS}, volume 33, 2020.

\bibitem{hybridminimax2020}
Quoc Tran-Dinh, Deyi Liu, and Lam M. Nguyen.
\href{https://papers.nips.cc/paper/2020/hash/7f141cf8e7136ce8701dc6636c2a6fe4-Abstract.html}{Hybrid Variance-Reduced SGD Algorithms for Minimax Problems with Nonconvex-Linear Function}.
In \emph{NeurIPS}, volume 33, 2020.

\bibitem{adgda2023}
Feihu Huang, Xidong Wu, and Zhengmian Hu.
\href{https://proceedings.mlr.press/v206/huang23a.html}{AdaGDA: Faster
adaptive gradient descent ascent methods for minimax optimization}.
In \emph{AISTATS}, PMLR 206:2365--2389, 2023.

\bibitem{dream2024}
Lesi Chen, Haishan Ye, and Luo Luo.
\href{https://proceedings.mlr.press/v238/chen24b.html}{An Efficient Stochastic Algorithm for Decentralized Nonconvex-Strongly-Concave Minimax Optimization}.
In \emph{AISTATS}, PMLR 238:1990--1998, 2024.

\bibitem{vi2022}
Ahmet Alacaoglu and Yura Malitsky.
\href{https://proceedings.mlr.press/v178/alacaoglu22a.html}{Stochastic
variance reduction for variational inequality methods}.
In \emph{Conference on Learning Theory (COLT)}, PMLR 178:778--816, 2022.

\end{thebibliography}
\end{document}